\documentclass[11pt,letterpaper]{amsart}
\usepackage[letterpaper,hmargin=1.25in,vmargin=1.5in]{geometry}
\usepackage{latexsym, amssymb, amsmath}
\usepackage{booktabs}   
\usepackage{multirow}
\usepackage{color}

\newtheorem{conjecture}{Conjecture}
\newtheorem{theorem}{Theorem}[section]
\newtheorem{lemma}[theorem]{Lemma}

\newtheorem{remark}[theorem]{Remark}
\theoremstyle{definition}

\makeatletter
    
    \@addtoreset{equation}{section}
  \makeatother

\begin{document}

\title[Classification of biharmonic submersions]
{Classification of biharmonic Riemannian submersions from manifolds with constant sectional curvature}

\author{Shun Maeta}
\address{Department of Mathematics, Chiba University, 1-33, Yayoicho, Inage, Chiba, 263-8522, Japan.}
\curraddr{}
\email{shun.maeta@faculty.gs.chiba-u.jp~{\em or}~shun.maeta@gmail.com}

\author{Miho Shito}
\address{Department of Mathematics, Chiba University, 1-33, Yayoicho, Inage, Chiba, 263-8522, Japan.}
\email{shito.miho@gmail.com}

\date{}

\dedicatory{}


\commby{}

\begin{abstract}
In 2011, Wang and Ou (Math. Z. {\bf 269}:917-925, 2011) showed that any biharmonic Riemannian submersion from a 3-dimensional Riemannian manifold with constant sectional curvature to a surface is harmonic. In this paper, we generalize the 3-dimensional setting to arbitrary dimensions. 
By constructing an adapted orthonormal frame, we simplify the biharmonic equation for Riemannian submersions and analyze the curvature properties of Riemannian manifolds with constant sectional curvature. 
As a result, we prove that a Riemannian submersion from an $(n+1)$-dimensional Riemannian manifold with constant sectional curvature to an $n$-dimensional Riemannian manifold is biharmonic if and only if it is harmonic. 
This result may also be viewed as an affirmative codimension-one Riemannian submersion analogue of Chen's conjecture, the generalized Chen's conjecture, and the BMO conjecture.

\end{abstract}

\maketitle
\bibliographystyle{amsalpha}

\section{Introduction}\label{intro}

In differential geometry, a biharmonic map is studied as a natural extension of a harmonic map. 
A harmonic map is defined as a critical point of the energy functional
\[
E(\phi) = \frac{1}{2} \int_M |d\phi|^2 \, d\mu_g,
\]
where $\phi: (M,g) \to (N,h)$ is a smooth map, and $d\mu_g$ is the volume element on $M$.

On the other hand, a biharmonic map is defined as a critical point of the bi-energy functional
\[
E_2(\phi) = \frac{1}{2} \int_M |\tau(\phi)|^2 \, d\mu_g,
\]
where $\tau(\phi)$ is the tension field of the map $\phi$, and a harmonic map satisfies $\tau(\phi) = 0$ (cf. \cite{EL2, ES}). 
Thus, if $\tau(\phi) = 0$, then $E_2(\phi) = 0$, so every harmonic map is a biharmonic map.

If an isometric immersion $\phi: (M,g) \to (N,h)$ is biharmonic, then $M$ is called a biharmonic submanifold in $N$.

The study of biharmonic maps includes several well-known open conjectures. 
We outline three such conjectures and summarize recent progress on each.

\begin{conjecture}[Chen's conjecture \cite{C}]
The only biharmonic submanifolds in Euclidean spaces are the minimal ones.
\end{conjecture}

\begin{remark}
This conjecture is customarily referred to as Chen's conjecture, rather than Chen conjecture, and we follow this convention in the present paper.
\end{remark}

Chen, Ishikawa, and Jiang independently proved that biharmonic surfaces in $\mathbb{E}^3$ are minimal \cite{C, CI, J2}. 
Hasanis, Vlachos, and Defever confirmed that biharmonic hypersurfaces in $\mathbb{E}^4$ are minimal \cite{D, HV}. 
Fu, Hong, and Zhan resolved the conjecture for hypersurfaces in $\mathbb{E}^5$ and $\mathbb{E}^6$ \cite{FZ1, FZ2}. 
Even for hypersurfaces in $\mathbb{E}^{n+1}$, the case $n \geq 6$ remains open.
Related works include \cite{AM, CM, Di, Di2, F, MOR, OC}.

\begin{conjecture}[Generalized Chen's conjecture \cite{CMO}]
Every biharmonic submanifold in a non-positively curved Riemannian manifold is minimal.
\end{conjecture}

There have been many studies on this conjecture \cite{BMO, BMO2, CMO, CMO2, FH, GV, J, Maeta14, NU, On, OC}, but Ou and Tang have constructed counterexamples \cite{OT}.
If the ambient space is a hyperbolic space $\mathbb{H}^{n+1}$, the conjecture is still open. 
There are many affirmative partial answers to this conjecture for surfaces and hypersurfaces in $\mathbb{H}^{n+1}$. 
Caddeo, Montaldo, and Oniciuc resolved the conjecture for surfaces in $\mathbb{H}^3$ \cite{CMO}.
Balmu\c{s}, Montaldo, and Oniciuc resolved the conjecture for hypersurfaces in $\mathbb{H}^4$ \cite{BMO2}. 
Guan, Li, and Vrancken proved the conjecture for hypersurfaces in $\mathbb{H}^{5}$ \cite{GV}. 
Fu, Hong, and Zhan proved the conjecture for hypersurfaces in $\mathbb{H}^{6}$ \cite{FZ2}.
However, even for hypersurfaces in $\mathbb{H}^{n+1}$, the case $n \geq 6$ remains open.

\begin{conjecture}[BMO conjecture \cite{BMO}]
Every biharmonic submanifold in a sphere has constant mean curvature.
\end{conjecture}

Caddeo, Montaldo, and Oniciuc proved the conjecture for surfaces in $\mathbb{S}^3$ \cite{CMO}.
Balmu\c{s}, Montaldo, and Oniciuc proved the conjecture for hypersurfaces in $\mathbb{S}^4$ \cite{BMO}.
Guan, Li, and Vrancken proved the conjecture for hypersurfaces in $\mathbb{S}^5$ \cite{GV}.
Fu, Hong, and Zhan proved the conjecture for hypersurfaces in $\mathbb{S}^6$ \cite{FZ2}.
However, even for hypersurfaces in $\mathbb{S}^{n+1}$, the case $n \geq 6$ remains open.
 Other related studies include \cite{BLO, FO, F2, LM17, Maeta17}.

Therefore, these conjectures remain unsolved even for hypersurfaces in Euclidean spaces $\mathbb{E}^{n+1}$, spheres $\mathbb{S}^{n+1}$, and hyperbolic spaces $\mathbb{H}^{n+1}$ for dimensions $n \geq 6$.

Riemannian submersions, considered the dual concept of isometric immersions, have been studied less extensively than their counterparts.

In 2011, Wang and Ou proved ~\cite{WO} the necessary and sufficient conditions for a Riemannian submersion from a 3-dimensional Riemannian manifold with constant sectional curvature to a surface to be a biharmonic map. 
Furthermore, they showed that any biharmonic Riemannian submersion from a 3-dimensional Riemannian manifold with constant sectional curvature to a surface is harmonic.
In their paper, the integrability data $\{f_1, f_2, \kappa_1, \kappa_2, \sigma\}$ were introduced to analyze the structure of Riemannian submersions (see Section 2 for details).
 The integrability data consist of functions on the manifold defined from the Lie brackets of a local orthonormal frame adapted to the Riemannian submersion.
  This framework allows for capturing information about connections and curvatures, simplifying the biharmonic equation
\[
\left( \nabla^{\phi}_{e_i} \nabla^{\phi}_{e_i} - \nabla^{\phi}_{\nabla^{M}_{e_i}e_i} \right)\tau(\phi)- R^N(d\phi(e_i), \tau(\phi))d\phi(e_i) = 0,
\]
which was derived by Jiang in 1986~\cite{J}.
Here, $\nabla^{\phi}$ denotes the induced connection, and $R^N$ denotes the Riemannian curvature tensor of $N$, defined by
\[
R^N(X, Y)Z = \nabla^N_X \nabla^N_Y Z - \nabla^N_Y \nabla^N_X Z - \nabla^N_{[X,Y]} Z,
\]
for vector fields $X, Y, Z$ on $N$, and $\{e_1, \dots, e_n\}$  is a local orthonormal frame field of $M$.

In 2019, Akyol and Ou generalized ~\cite{AO} the integrability data to the form
\[
\{f_{ij}^k, \kappa_i, \sigma_{ij}\}_{i,j,k = 1,\dots,n}
\]
and computed the biharmonic equation for Riemannian submersions in general dimensions as follows (detailed in Section~2):
\begin{equation}\label{2019}
\begin{aligned}
&\sum_{i=1}^{n+1} e_i e_i(\kappa_k) - \sum_{i,j=1}^{n} P_{ii}^j e_j(\kappa_k) - \sum_{i=1}^{n} \kappa_i e_i(\kappa_k) \\
&\quad + \sum_{i,j=1}^n \left[ 2 e_i(\kappa_j) P_{ij}^k + \kappa_j (e_i P_{ij}^k) + \kappa_j P_{ij}^l P_{il}^k - \kappa_i \kappa_j P_{ij}^k - \kappa_j P_{ii}^l P_{lj}^k \right] \\
&\quad + \text{Ricci}^N(d\phi(\mu), d\phi(e_k)) = 0, \quad k = 1, 2, \dots, n,
\end{aligned}
\end{equation}
where $P^k_{ij} = \frac{1}{2} \left( - f^j_{ik} - f^i_{jk} + f^k_{ij} \right)$ for all $i, j, k = 1, 2, \dots, n$.

Wang and Ou have recently conducted a series of studies~\cite{WO2, WO3, WO4, WO5} that classify proper biharmonic Riemannian submersions from 3-dimensional manifolds, including Thurston's eight model geometries, Bianchi-Cartan-Vranceanu (BCV) spaces, and product manifolds such as $M^2 \times \mathbb{R}$ (see \cite{Ou} for details).
 However, the classification in general dimensions remains largely unexplored.

In this paper, we establish the necessary and sufficient conditions for a biharmonic Riemannian submersion from an $(n+1)$-dimensional  Riemannian manifold with constant sectional curvature to be a harmonic map.
The main theorem is as follows.

\begin{theorem}\label{main}
For $n\geq2$, let 
$
\phi : (M^{n+1}(c), g) \longrightarrow (N^n, h)
$
be a Riemannian submersion from a Riemannian manifold with constant sectional curvature c. 
Then $ \phi $  is biharmonic if and only if it is harmonic.
\end{theorem}
This result may also be viewed as an affirmative codimension-one Riemannian submersion analogue of Chen's conjecture, the generalized Chen's conjecture, and the BMO conjecture.

\begin{remark}
It is known that the map $\phi:(\mathbb{H}^2,\frac{dx^2+dy^2}{y^2})\rightarrow ((0,\infty),\frac{dy^2}{y^2}), \phi(x,y)=y$ is a proper biharmonic Riemannian submersion (see, for example, \cite{Usula26}). Thus, the assumption $n\geq2$ is sharp.
\end{remark}

In order to prove Theorem \ref{main}, we needed to overcome several difficulties. 
In the following, we describe how these challenges were addressed and outline the steps of the proof.
\\
\\
\textbf{Step 1: Computation of $ R^M_{abcd} $ and $ \nabla_{e_a} e_b $. 
}

In the 3-dimensional case ($n=2$), the integrability data consist of 5 components: 
\[
\{f_1, f_2, \kappa_1, \kappa_2, \sigma\}.
\] 
Although this complicates the biharmonic equation for a Riemannian submersion, Wang and Ou \cite{WO}  succeeded in their analysis by considering the components of the curvature tensor $R^M_{1312},$ $R^M_{1313},$ $R^M_{1323},$ $R^M_{1212},$ $R^M_{1223},$ $R^M_{2313},$ and $R^M_{2323}$, and by exploiting the fact that the integrability data consist of only five components.

However, even in the 4-dimensional case, the integrability data increase to 15 components  
\[
\{\kappa_1, \kappa_2, \kappa_3, f_{12}^1, f_{12}^2, f_{12}^3, f_{13}^1, f_{13}^2, f_{13}^3, f_{23}^1, f_{23}^2, f_{23}^3, \sigma_{12}, \sigma_{13}, \sigma_{23}\},
\]
the connection terms $\nabla_{e_i} e_j$ require computing 16 combinations, $P_{ij}^k$ require computing 27 combinations, and the computation of the curvature tensor $R^M_{abcd}$ requires at least 21 terms.

We remark that in $n+1$ dimensions, the number of components of the integrability data increase to as many as $n + \frac{n^2(n-1)}{2} + \frac{n(n-1)}{2}$.

To overcome this difficulty, we succeeded in identifying the following four components of the curvature tensor and the following four connection terms as the only ones required.
\[
 R^M_{a (n+1) c d} ,  \quad R^M_{a b a b} , \quad R^M_{a (n+1) a (n+1)} , \quad R^M_{a (n+1) c (n+1)}
\]
and
 \[
  \nabla_{e_a} e_b , \quad \nabla_{e_a} e_{n+1} , \quad \nabla_{e_{n+1}} e_a , \quad \nabla_{e_{n+1}} e_{n+1} ,
  \]
 where $ a, b, c, d = 1, \dots, n $.\\

\noindent
\textbf{Step 2: Proof of constancy of integrability data along $ e_{n+1} $.}

To analyze the biharmonic equation for a Riemannian submersion $\phi$, we need to show that for an orthonormal frame $\{e_1,\cdots,e_n, e_{n+1}\}$ adapted to a Riemannian submersion $\phi$ with $e_{n+1}$ being vertical, all the integrability data are constant along fibers of $\phi$, that is, $ e_{n+1}(f^c_{ab}) = e_{n+1}(\kappa_a) = e_{n+1}(\sigma_{ab}) = 0 $. 
In the 3-dimensional case, Wang and Ou showed it (cf. Lemma 3.1 in \cite{WO}). 
Unfortunately, even in the 4-dimensional case, one cannot show it in general. 
However, interestingly, the biharmonicity condition overcomes this difficulty, that is, if $\phi$ is biharmonic, then we can show that $ e_{n+1}(f^c_{ab}) = e_{n+1}(\kappa_a) = e_{n+1}(\sigma_{ab}) = 0 $ (Lemma \ref{lemma3.1}). 
In fact, using the definitions and results from Step 1, the equalities $ e_{n+1}(f^c_{ab}) = e_{n+1}(\sigma_{ab}) = 0 $ can be readily established. For the more challenging proof of $ e_{n+1}(\kappa_a) = 0 $, we employ the biharmonic equation \eqref{2019}, simplified as
\begin{equation*}
e_{n+1} e_{n+1} (\kappa_a) = -\kappa_a \left( \sum_{i=1}^n \kappa_i^2 + (2n-1)c \right),
\end{equation*}
and combine this with previous results to complete the proof of Lemma \ref{lemma3.1}.
This result simplifies computations and allows us to derive additional usable expressions from the curvature tensor $ R^M_{abcd} $. \\

\noindent
\textbf{Step 3: Construction of an adapted orthonormal frame (Lemma \ref{lemma3.2}).}

Although Lemma \ref{lemma3.1} shows the constancy of the integrability data along $e_{n+1}$, this is still far from sufficient to conclude that a biharmonic Riemannian submersion is harmonic. In fact, the equation for biharmonic Riemannian submersions is given by (\ref{2019}).

To overcome this difficulty, we construct a special orthonormal frame satisfying
\begin{equation*}
\kappa_2 = \kappa_3 = \dots = \kappa_n = 0
\end{equation*}
and, for any $ i = 1, \dots, n-2 $,
\begin{equation*}
\sigma_{i,i+j} = 0, \quad (j \geq 2).
\end{equation*}

To construct this frame, it was necessary to maintain orthonormality while satisfying the conditions required by the integrability data, specifically, the conditions on the Lie bracket of an orthonormal frame compatible with a Riemannian submersion. 
The Gram-Schmidt orthogonalization method proved insufficiently precise for this purpose, so we successfully constructed the frame by repeatedly applying Householder transformations.\\

\noindent
\textbf{Step 4: Completion of the proof of Theorem \ref{main} by contradiction.}

Lemma \ref{lemma3.2} simplifies the biharmonic equation as follows.
\begin{equation*}
\kappa_1^2 + \sigma_{12}^2 - 4 \sum_{i=1}^{n-1} \sigma_{i,i+1}^2 + (2n-1)c = 0.
\end{equation*}
We aim to show that $ \kappa_1 = 0$. 
To show that $\kappa_1=0$, an analysis of $\sigma_{i,i+1}$ is required.
We overcome this difficulty by differentiating along all directions of the adapted frame and ultimately show that a biharmonic Riemannian submersion is harmonic.

\section{Preliminaries}
In this section, we recall some definitions and facts about biharmonic Riemannian submersions which will be used in this paper. 

Let $\phi : (M^{n+1}, g) \to (N^n, h)$ be a Riemannian submersion.  
A local orthonormal frame adapted to $\phi$ consists of horizontal lifts of vector fields from the base manifold $N$, which locally span the horizontal distribution on $M$.  
Since basic vector fields locally span the horizontal distribution, such a frame can always be constructed (cf. page 2 in \cite{Mari}).

Let $\{e_1, \dots, e_n, e_{n+1}\}$ be an orthonormal frame adapted to the submersion $\phi$, where $e_1, \dots, e_n$ are horizontal lifts of the orthonormal frame $\{\varepsilon_1, \dots, \varepsilon_n\}$ on $N$, and $e_{n+1}$ is vertical along the fiber.

The Lie bracket $[e_i, e_{n+1}]$ is vertical and does not contribute to the horizontal component 
(Lemma 3, \cite{ON}).
Moreover, the bracket $[e_i, e_j]$ of horizontal lifts corresponds under $\phi$ to the bracket $[\varepsilon_i, \varepsilon_j]$ on $N$, preserving the Lie bracket structure
(Lemma 1, \cite{ON}).

In \cite{WO}, the integrability data were introduced.
Subsequently, in \cite{AO}, the integrability data for the adapted frame of a Riemannian submersion $\phi$ and the biharmonic equation were generalized as follows:\\
If we assume that
\begin{equation}
[\varepsilon_i, \varepsilon_j] = F^k_{ij} \varepsilon_k, 
\end{equation}
where, in the sequel, $F^k_{ij} \in C^\infty(N)$ and the Einstein convention is used, then we have
\begin{equation}\label{data}
\begin{aligned}
&[e_i, e_{n+1}] = \kappa_i e_{n+1},  \\
&[e_i, e_j] = f^k_{ij} e_k - 2 \sigma_{ij} e_{n+1}, \quad i, j = 1, 2, \dots, n, 
\end{aligned}
\end{equation}
where $f^k_{ij} = F^k_{ij} \circ \phi$, $\kappa_i$, and  $\sigma_{ij} \in C^\infty(M)$ for all  $i, j = 1, 2, \dots, n$. We will call $\{f^k_{ij}, \kappa_i, \sigma_{ij}\}$ the  integrability data of the adapted frame of the Riemannian submersion $\phi$. It follows from (\ref{data}) that
\begin{equation}
f^k_{ij} = - f^k_{ji}, \quad \sigma_{ij} = - \sigma_{ji}, \quad ( i, j, k = 1, 2, \dots, n). 
\end{equation}

\begin{theorem}[Theorem 3.1 in \cite{AO}]\label{theoremAO}
Let $\phi : (M^{n+1}, g) \to (N^n, h)$ be a Riemannian submersion with the adapted frame $\{e_1, \dots, e_{n+1}\}$ and the integrability data  $\{f^k_{ij}, \kappa_i, \sigma_{ij}\}$.
Then, the Riemannian submersion $\phi$ is biharmonic if and only if
\begin{equation}
\begin{aligned}
&\sum_{i=1}^{n+1} e_i e_i(\kappa_k) - \sum_{i,j=1}^n P_{ii}^j e_j(\kappa_k) - \sum_{i=1}^n \kappa_i e_i(\kappa_k)\\
&\quad  + \sum_{i,j=1}^n \left[ 2 e_i(\kappa_j) P_{ij}^k + \kappa_j (e_i P_{ij}^k) + \kappa_j P_{ij}^l P_{il}^k - \kappa_i \kappa_j P_{ij}^k - \kappa_j P_{ii}^l P_{lj}^k \right]  \\
&\quad  + \text{Ricci}^N(d\phi(\mu), d\phi(e_k)) = 0,\quad k = 1, 2, \dots, n,
\end{aligned}
\label{eq:biha}
\end{equation}
where 
\[
P^k_{ij} = \frac{1}{2} \left( - f^j_{ik} - f^i_{jk} + f^k_{ij} \right) \text{ for all }  i, j, k = 1, 2, \dots, n,
\]
and 
\[
\mu=(\nabla^M_{e_{n+1}} e_{n+1})^\mathcal{H}.
\]
\begin{proof}
See the proof of Theorem 3.1 in \cite{AO}.
\end{proof}
\end{theorem}
This equation  (\ref{eq:biha}) is called the biharmonic equation.

In this paper, we use both $\sigma_{ij}$ and $\sigma_{i,j}$ to denote the same function defined on a manifold. 
The comma is used in expressions like $\sigma_{i,j+1}$ for clarity, especially when indices involve arithmetic operations.

\section{Proof of the Main Theorem and Supporting Lemmas}
In this section, we present three lemmas that play essential roles in the proof of Theorem~\ref{main}, and show Theorem~\ref{main}.  
Each lemma focuses on a different aspect of the geometric or analytic structure involved in the main result.

We begin by showing that the integrability data remain constant along the fiber direction $ e_{n+1} $.  
This property will be used directly in the proof of Theorem~\ref{main}.
\begin{lemma}\label{lemma3.1}
Let $\phi: (M^{n+1}(c),g) \to (N^n, h)$ be a biharmonic Riemannian submersion from an $(n+1)$-dimensional Riemannian manifold $M^{n+1}(c)$ of constant sectional curvature $c$ to an arbitrary $n$-dimensional Riemannian manifold $(N^n, h)$. For any orthonormal frame $\{e_1, \dots, e_n, e_{n+1}\}$ on $M^{n+1}(c)$ adapted to the Riemannian submersion $\phi$, with $e_{n+1}$ being a vertical vector field, all integrability data $f^k_{ij}$, $\kappa_i$, and $\sigma_{ij}$ are constant along the fibers of $\phi$, that is, the following equation holds for all $i,j,k = 1, 2, \dots, n$:
\[
e_{n+1}(f^k_{ij}) = e_{n+1}(\kappa_i) = e_{n+1}(\sigma_{ij}) = 0.
\]
\begin{proof} 
Since the case for $n=2$ has already been proven (cf. \cite{WO}), we will show the case for $n \geq 3$.
Since $e_{n+1}$ is a vertical vector field, we obtain $d\phi(e_{n+1}) = 0$. \\
Hence, by differentiating the function $f^k_{ij} = F^k_{ij} \circ \phi$ by $e_{n+1}$ we get, 
\begin{align*}
 e_{n+1}(f^k_{ij}) &= e_{n+1}( F^k_{ij} \circ \phi) \\
 &= (d\phi(e_{n+1}))( F^k_{ij} )\\
 &= 0.
\end{align*}

Therefore, for any $i,j,k = 1, 2, \dots, n$, we have
\begin{align*}
 e_{n+1}(f^k_{ij}) = 0, \quad e_{n+1}(P_{ij}^k) = 0. \label{eq:result}
\end{align*}

Next, we show that $e_{n+1}(\sigma_{ij}) = 0$. 
A straightforward computation using (\ref{data}) and Koszul formula gives

\begin{equation}\label{nabla}
\begin{aligned}
\nabla^M_{e_i} e_j &= P^k_{ij} e_k - \sigma_{ij} e_{n+1} \quad  \text{for any }  i,j,k = 1, 2, \dots, n,\\
\nabla^M_{e_{n+1}} e_{n+1} &= \sum_{i=1}^n \kappa_i e_i,\\
\nabla^M_ { e _ { i } } e _ { n + 1 } &= \sigma _ { i j } e _ { j }\quad \text{for any }  i = 1, 2, \dots, n,\\
\nabla^M_ { e _ { n+1 } } e _ { i } &=\sigma _ { i j } e _ { j } -\kappa_i e_{n+1}\quad  \text{for any }  i = 1, 2, \dots, n.
\end{aligned}
\end{equation}
Using  (\ref{nabla}) and the fact that $M^{n+1}(c)$ has constant curvature, we get 
\begin{equation}
\left\{ \,
\begin{aligned}
&- R^M _ { a (n + 1)c d} = e _ { a } ( \sigma _ { c d } ) +P_ { a  l } ^ { d } \sigma _ { c l}  - P_ { a c } ^ { l} \sigma _ { ld }- \kappa _ { c } \sigma _ { a d }  + \kappa _ { d } \sigma _ { a c } - \kappa _ { a } \sigma _ { cd } = 0, \\
&- R^M _ { a b a b } = e _ { a } ( P_ { b a } ^ { b } ) + P_ { b a } ^ { l} P_ { a l } ^ { b } + 3 \sigma _ { a b } ^ { 2 } - e _ { b } ( P _ { a a } ^ { b } ) -P_{aa}^l P_{bl}^b- f_ { a b } ^ { l } P_ { la } ^ { b }=-c, \\
&- R^M _ { a (n + 1 ) a (n + 1 ) } = - \sigma _ { a l } ^ { 2 } - e _ { a } ( \kappa _ { a } ) + P_ { a a } ^ { l } \kappa _ { l } + \kappa _ { a } ^ { 2 } = - c, \\
&- R^M _ { a (n + 1 )c (n + 1) } = - \sigma _ { c l } \sigma _ { a l } - e _ { a  } (\kappa _ { c } )+ P_ { a c } ^l\kappa_l+ \kappa _ { a }\kappa_c + e _ { n + 1 } ( \sigma _ { a c } ) = 0 \quad (a\not=c). \\ 
\end{aligned}
\label{eq:R^M}
\right.
\end{equation}

By differentiating both sides of the second equation in (\ref{eq:R^M}) by $e_{n+1}$, and substituting
\begin{align*}
e_{ n + 1 }  e_a = -[e_a, e_{n+1}] +e_a e_{ n + 1 } =-\kappa_a e_{ n + 1 }  + e_a e_{ n + 1 },
\end{align*}
we obtain 
\begin{align*}
 \sigma _ { a b } e _ { n + 1 } ( \sigma _ { a b } ) = 0.          
\end{align*}

Let
\[
\Omega_1=\{p\in M^{n+1}(c): e_{n+1}(\sigma_{ab})(p)\neq0\}.
\]
If $\Omega_1$ is nonempty, then it is open. 
We arrive at a contradiction.
Under this assumption, we have 
\[
\sigma_ { ab} = 0
\]
 on $\Omega_1$.
This contradicts the assumption that  $ e_{n+1}(\sigma_{ab}) \neq 0 $.
Therefore, we conclude that \( e_{n+1}(\sigma_{ab}) = 0 \). 

Next, we show  $e_{n+1}(\kappa_a) = 0$. 
As before, by differentiating both sides of the first equation in (\ref{eq:R^M}) by $e_{n+1}$, and substituting 
\(
e_{ n + 1 }  e_a = -\kappa_a e_{ n + 1 }  + e_a e_{ n + 1 } 
\)
we obtain
\begin{align*}
- \sigma _ { a d } e _ { n + 1 } ( \kappa _ { c } ) + \sigma _ { a c } e _ { n + 1 } ( \kappa _ { d } ) - \sigma _ { c d } e _ { n + 1 } ( \kappa _ { a } ) = 0 \quad \text{for any }a,c,d= 1, 2, \dots, n .
\end{align*}

First, we relabel the indices by replacing $ a $ with $ i $, $ c $ with $ j $, and $ d $ with $ k $. Under this assignment, the following identity holds
\begin{equation}\label{eq:ijk}
\sigma_{ij} e_{n+1}(\kappa_k) = \sigma_{ik} e_{n+1}(\kappa_j) + \sigma_{jk} e_{n+1}(\kappa_i).
\end{equation}

Next, we relabel the indices by replacing $a$ with $j$, $c$ with $k$, and $d$ with $i$.
Since $\sigma = (\sigma_{ij})$ is skew-symmetric, the following identity holds
\begin{equation}\label{eq:jki}
\sigma_{ij} e_{n+1}(\kappa_k) = -\sigma_{ik} e_{n+1}(\kappa_j) - \sigma_{jk} e_{n+1}(\kappa_i).
\end{equation}
By combining equations (\ref{eq:ijk}) with (\ref{eq:jki}), we obtain
\[
\sigma _ { i j } e _ { n+1 }(\kappa_k)=0
\]  
for any $ i,j,k = 1, 2, \dots, n$. 

Assume that there exists a fixed index $a$ and a nonempty open set
$\Omega_2 \subset M^{n+1}(c)$ such that
\[
e_{n+1}(\kappa_a)\neq0
\]
on $\Omega_2$.
Then we arrive at a contradiction.
Under this assumption, we have 
\[
\sigma _ { i j }= 0
\]
 on $\Omega_2$.
 Under this condition, equation  (\ref{eq:R^M}) becomes
\begin{equation}
\left\{
\begin{aligned}
&e _ { a } ( P_ { b a } ^ { b } ) - e _ { b } ( P_ { a a } ^ { b } ) = - c + P_ { a a } ^ { l } P _ { b l } ^ { b } + f _ { a b } ^ { l } P _ { la } ^ { b } - P_ { b a } ^ { l} P_ { a l } ^ { b },\\
&e _ { a } ( \kappa _ { c } ) = \delta _ { a c } c + P_ { a c } ^ { l } \kappa _ { l } + \kappa _ { a } \kappa_ { c }.\\
\end{aligned}
\label{eq:R^Msigma=0}
\right.
\end{equation}

By substituting the second equation of (\ref{eq:R^Msigma=0}) into (\ref{eq:biha}), each term becomes as follows.
\begin{equation}
\begin{aligned}
&\quad\sum _ { i = 1 } ^ { n + 1 } e _ { i } e _ { i } ( \kappa _ { k } )=\sum _ { i = 1 } ^ { n } \{ P _ { i k } ^ { l } P _ { i l } ^ { m } \kappa _ { m } + P _ { i k } ^ { l }\kappa_i\kappa_l+ \kappa_l e _ { i } ( P _ { i k } ^ { l } )+ P _ { i k } ^ { l }\kappa_i\kappa_l + P _ { i i } ^ { l }\kappa_k\kappa_l +2\kappa_i^2\kappa_k\}\\
&\quad \qquad \qquad\qquad\qquad \qquad+ ( n + 1 ) \kappa _ { k } c + P _ { l k } ^ { l } c + e _ { n + 1 } e _ { n + 1 } ( \kappa _ { k } ), \\
&- \sum _ { i , j = 1 } ^ { n } P _ { i i } ^ { j } e _ { j } ( \kappa _ { k } ) =  \sum _ { i , j = 1 } ^ { n }\{ - P _ { i i } ^ { j } P _ { j k } ^ { l } \kappa _ { l } - P _ { i i } ^ { j } \kappa _ { j } \kappa_ { k } \}+ \sum _ { i = 1 } ^ { n } \{ - P _ { i i } ^ { k } c \},\\
&- \sum _ { i = 1 } ^ { n } \kappa _ { i } e _ { i } ( \kappa _ { k } ) = \sum _ { i = 1 } ^ { n } \{ - P _ { i k } ^ { l } \kappa _ { i } \kappa _ { l } - \kappa _ { i } ^ { 2 } \kappa _ { k } \} - \kappa _ { k } c,\\
&\sum_{i,j=1}^n \{ 2 e_i(\kappa_j) P_{ij}^k + \kappa_j e_i( P_{ij}^k) + \kappa_j P_{ij}^l P_{il}^k - \kappa_i \kappa_j P_{ij}^k - \kappa_j P_{ii}^l P_{lj}^k \} \\
&\qquad =\sum_{i,j=1}^n \left\{2P_{ij}^kP_{ij}^l\kappa_l +2P_{ij}^k\kappa_i \kappa_j+\kappa_j e_i (P_{ij}^k)\right.\\
&\qquad\qquad\qquad\qquad 
\left.+ P_{ij}^l P_{il}^k\kappa_j  - P_{ij}^k\kappa_i \kappa_j  -P_{ii}^l P_{lj}^k \kappa_j  \right\}+\sum_{i=1}^n2P_{ii}^k c, \\
&  \mathrm{Ricci}^ { N } ( d \phi ( \mu ) , d \phi ( e _ { k } ) ) =(n-1)\kappa_kc.
\end{aligned}
\end{equation}
By simplifying these expressions using the fact that $P _ { i j } ^ { k } = - P _ { i k } ^ { j }$, we obtain the following equation:
\begin{equation}
e _ { n + 1 } e _ { n + 1 } ( \kappa _a ) = - \kappa _ { a } (  \sum _ { i = 1 } ^ { n }\kappa _ { i } ^ { 2 } + ( 2 n - 1 ) c ). \label{eq:eek}
\end{equation}\\
To derive additional usable equations, we differentiate both sides of the second and third equations in (\ref{eq:R^M}) by $e_{ n + 1 }e_{ n + 1 }$.
Then, by substituting $\sigma_{ij} = 0$ and $e_{ n + 1 }  e_a = -\kappa_a e_{ n + 1 }  + e_a e_{ n + 1 } $ we obtain 
\begin{equation}
\left\{ \,
\begin{aligned}
&- 3 \{ e _ { n + 1 } ( \kappa _ { a } ) \} ^ { 2 } - 4 \kappa _a e _ { n + 1 } e _ { n + 1 } ( \kappa _a ) - P_ { a a } ^ { l } e _ { n + 1 } e _ { n + 1 }( \kappa _l ) + e _ { a } e _ { n + 1 } e _ { n + 1 } ( \kappa _a ) = 0, \\
&- 3e _ { n + 1 } ( \kappa_a ) e _ { n + 1 } ( \kappa _ { c} ) - 3 \kappa _a e _ { n + 1 } e _ { n + 1 } ( \kappa _ { c} ) - \kappa _c e _ { n + 1 } e _ { n + 1 } ( \kappa _ { a} )\\
&\qquad \qquad\qquad\qquad\qquad -P_{ac}^le _ { n + 1 } e _ { n + 1 }( \kappa _ { l} ) + e _ { a } e _ { n + 1 } e _ { n + 1 } ( \kappa _ { c } ) = 0.\\
\end{aligned}
\label{eq:eeR^Msigma=0}
\right.
\end{equation}
By substituting (\ref{eq:R^Msigma=0}) and (\ref{eq:eek}) into (\ref{eq:eeR^Msigma=0}), we obtain the following equation:
\begin{equation}
\left\{ \,
\begin{aligned}
&- 3 \{ e _ { n + 1 } ( \kappa _ { a } ) \} ^ { 2 } + \kappa _ { a } ^ { 2 } L - K c=0, \\
&- 3 \{ e _ { n + 1 } ( \kappa _ { b } ) \} ^ { 2 } + \kappa _ { b } ^ { 2 } L - K c=0, \\
&3e _ { n + 1 } ( \kappa _ { a } ) e _ { n + 1 } ( \kappa _ { b } ) = \kappa _ { a } \kappa _ { b } L,
\end{aligned}
\label{eq:eek^2}
\right.
\end{equation} 
where
\[
L= \sum _ { i = 1 } ^ { n }\kappa_i^2+(6n-5)c,~~K= \sum _ { i = 1 } ^ { n }\kappa _ { i } ^ { 2 }  + ( 2 n - 1 ) c.
\]

\noindent
Case 1.  $c=0$: In this case, from  (\ref{eq:eek}) and (\ref{eq:eek^2}), we have
\begin{equation}\label{c=0}
\begin{aligned}
3 \{ e _ { n+1 } ( \kappa _ { a } ) \} ^ { 2 } &= \sum _ { i = 1 } ^ { n }\kappa_a^2\kappa_i^2  ,\\
e _ { n + 1 } e _ { n + 1 } ( \kappa_a ) &= - \sum _ { i = 1 } ^ { n }\kappa_a\kappa_i^2 , \\
 e _ { n + 1 } ( \kappa_b )&= \sum _ { i = 1 } ^ { n }\frac{\kappa_a \kappa_b \kappa_i^2}{3e_{n+1}(\kappa_a)}.
 \end{aligned}
\end{equation}
Differentiating both sides of the first equation in \eqref{eq:eek^2} by $e _ { n + 1 }$, we obtain
\[
3e_{n+1}(\kappa_a)e_{n+1}e_{n+1}(\kappa_a)- \sum _ { i = 1 } ^ { n }\kappa_a\kappa_i^2e_{n+1}(\kappa_a)- \sum _ { i = 1 } ^ { n }\kappa_a^2\kappa_ie_{n+1}(\kappa_i)=0.
\]
Substituting the second equation of \eqref{c=0}  into this, we obtain
\[
4\sum _ { i = 1 } ^ { n }\kappa_a\kappa_i^2e_{n+1}(\kappa_a)+\sum _ { i = 1 } ^ { n }\kappa_a^2\kappa_ie_{n+1}(\kappa_i)=0.
\]
Since $\kappa_a\neq0$, it follows that 
\[
4\sum _ { i = 1 } ^ { n }\kappa_i^2e_{n+1}(\kappa_a)+\sum _ { i = 1 } ^ { n }\kappa_a\kappa_ie_{n+1}(\kappa_i)=0.
\]
Substituting the first and third equations of \eqref{c=0}, we obtain
\[
4\sum _ { i = 1 } ^ { n }\kappa_i^2e_{n+1}(\kappa_a)+\sum _ { i,k = 1 } ^ { n }\kappa_a\kappa_i\frac{\kappa_a \kappa_i \kappa_k^2}{3e_{n+1}(\kappa_a)}=0.
\]
By simplifying this, we obtain 
\[
 12\{ e _ { n+1 } ( \kappa _ { a } ) \} ^ { 2 } +\sum_{k=1}^n\kappa_a^2 \kappa_k^2=0.
\]
Therefore, we obtain $e _ { n+1 } ( \kappa _ { a }) =0$, which contradicts the assumption.\\
\noindent
Case 2. $c\not=0$: 
Taking the difference between the first and second equations in (\ref{eq:eek^2}), we obtain
\begin{equation}
\left\{ \,
\begin{aligned}
&( \kappa _ { a } ^ { 2 } - \kappa _ { b } ^ { 2 } ) L - 3 ( \{ e _ { n+1 } ( \kappa _ { a } ) \} ^ { 2 } - \{ e _ { n+1 } ( \kappa _ { b } ) \} ^ { 2 } ) = 0, \\
&\kappa _ { a } \kappa _ { b } L - 3e _ { n+1 } ( \kappa _ { a } ) e _ { n+1 } ( \kappa _ { b } ) = 0.
\end{aligned}
\label{eq:eek^22}
\right.
\end{equation}

Since $e_{n+1}(\kappa_a) \neq 0$, we solve the second equation of (\ref{eq:eek^22}) for $e_{n+1}(\kappa_b)$, and substitute it into the first equation of (\ref{eq:eek^22}), which yields 
\[
\{ 3\{e _ { n+1 } ( \kappa _ { a } ) \}^2-\kappa _ { a} ^ { 2 } L \} \{ 3\{e _ { n+1 } ( \kappa _ { a } )\}^2 +\kappa _ { b } ^ { 2 } L \} = 0.
\]

Assuming that 
\[
3\{e _ { n+1 } ( \kappa _ { a } )\}^2 +\kappa _ { b } ^ { 2 } L
\not=0
\]
hold on 
$\Omega_3\subset \Omega_2$, it follows that 
\[
3 \{ e _ { n+1 } ( \kappa _ { a } ) \} ^ { 2 } - \kappa _ { a } ^ { 2 } L=0
\]
on $\Omega_3$. Substituting this into \eqref{eq:eek^2}, we obtain  
\[
( \sum _ { i = 1 } ^ { n }\kappa _ { i } ^ { 2 }  + ( 2 n - 1 ) c ) c =0.
\]
Since $c \neq 0$, we have
 \[
\sum _ { i = 1 } ^ { n } \kappa _ { i } ^ { 2 }  + ( 2 n - 1 ) c  =0.
 \]
Substituting this into (\ref{eq:eek}) and (\ref{eq:eek^2}), we obtain
\[
3 \{ e _ { n+1 } ( \kappa _ { a } ) \} ^ { 2 } = \kappa _ { a } ^ { 2 } 4(n-1)c ,\quad e _ { n + 1 } e _ { n + 1 } ( \kappa _a ) = 0.
\]
Differentiating the first equation by $e _ { n + 1 }$ and solving it together with the second equation, we obtain  $e _ { n+1 } ( \kappa _ { a }) =0$, which contradicts the assumption.
Therefore, 
\begin{equation}\label{case2ab}
3\{e _ { n+1 } ( \kappa _ { a } )\}^2 +\kappa _ { b } ^ { 2 } L=0
\end{equation}
on $\Omega_2$.
Moreover, since $b(\not=a)$  is arbitrary, it follows that 
\begin{equation}\label{case2ac}
3\{e _ { n+1 } ( \kappa _ { a } )\}^2 +\kappa _ { c } ^ { 2 } L=0
\end{equation}
holds for any $c\not=b,(c\not=a)$.
Solving \eqref{case2ab} and  \eqref{case2ac} simultaneously, we obtain 
\[
\kappa_b^2=\kappa_c^2
\]
on $\Omega_2$. 

When $c>0$, the equation \eqref{case2ab}  together with $L>0$ leads to a contradiction.

We consider the case where $c<0$.
Assume that there exists a subdomain $\Omega_3\subset\Omega_2$ such that $L<0$ holds within $\Omega_3$.
Squaring the third equation of \eqref{eq:eek^2} and substituting (\ref{case2ac}), we obtain
\[
\kappa_b^2L(3(e_{n+1}(\kappa_b))^2+\kappa_a^2L)=0
\]
on $\Omega_3$.
Assuming 
\[
3(e_{n+1}(\kappa_b))^2+\kappa_a^2L\not=0
\]
holds on $\Omega_4\subset\Omega_3$, it follows that 
\[
\kappa_b^2L=0
\]
 on $\Omega_4$. Substituting this into \eqref{case2ab}, we obtain $e_{n+1}(\kappa_a)=0$ on $\Omega_4$, which leads to a contradiction.
 Therefore,
\begin{equation}\label{case2ba}
3(e_{n+1}(\kappa_b))^2+\kappa_a^2L=0
\end{equation}
on $\Omega_3$. 

Assume that $\kappa_a\not=0$ holds on some subdomain $\Omega_4\subset\Omega_3$.
Then $e_{n+1}(\kappa_b)\not=0$ on $\Omega_4$, and by the same argument as above, we obtain 
\[
\{ 3\{e _ { n+1 } ( \kappa _ { b } ) \}^2-\kappa _ { b} ^ { 2 } L \} 
\{ 3\{e _ { n+1 } ( \kappa _ { b } )\}^2 +\kappa _ { c } ^ { 2 } L \} = 0,
\]
on $\Omega_4$ for $c\not=b$.
Assume that $3\{e _ { n+1 } ( \kappa _ { b } ) \}^2-\kappa _ { c } ^ { 2 } L\neq0$ on $\Omega_5\subset\Omega_4$.
Then, by the same argument as before, we arrive at a contradiction.
Therefore, 
\[
3\{e _ { n+1 } ( \kappa _ { b } )\}^2 =-\kappa _ { c } ^ { 2 } L=-\kappa _ { b } ^ { 2 } L
\]
on $\Omega_4$.
Combining this with \eqref{case2ba}, we obtain 
\[
\kappa_a^2=\kappa_b^2
\]
on $\Omega_4$.
Differentiating both sides of this equation by $e_{n+1}$, we obtain 
\[
\kappa_ae_{n+1}(\kappa_a)=
\kappa_be_{n+1}(\kappa_b),
\]
on $\Omega_4$.
Moreover, we have 
\[
\kappa_1^2=\kappa_2^2=\cdots=\kappa_n^2\not=0
\]
on $\Omega_4$.
Then, by multiplying both sides of the third equation in \eqref{eq:eek^2} by $\kappa_a\kappa_b$, we obtain
\[
3\kappa_a e_{n+1}(\kappa_a)\kappa_be_{n+1}(\kappa_b)=\kappa_a^2\kappa_b^2L.
\]
Therefore, we obtain 
\[
\frac{3}{4}
\{e_{n+1}(\kappa_a^2)\}^2=\kappa_a^4L
\]
which leads to a contradiction.

Therefore, $\kappa_a=0$ on $\Omega_3$. However, this contradicts the assumption that 
$e_{n+1}(\kappa_a)\not=0$ on $\Omega_2(\supset \Omega_3)$.

That is, the region where $L<0$ does not exist on $\Omega_2$.
Therefore, $L\geq0$ on  $\Omega_2$. 
However, in this case, it follows from \eqref{case2ab} that $e_{n+1}(\kappa_a)=0$ on $\Omega_2$, which leads to a contradiction.

Since a contradiction arises in every case, the assumption that $e _ { n+1 } ( \kappa _ { a }) \neq 0$ must be incorrect.

Therefore, $e _ { n+1 } ( \kappa _ { a }) =0$,
which implies that for any $i= 1, 2, \dots, n$, we have $e _ { n+1 } ( \kappa _ { i }) =0$.\\
Thus, the lemma is proved, and we obtain
\begin{equation}
e_{n+1}(f^k_{ij}) = e_{n+1}(\kappa_i) = e_{n+1}(\sigma_{ij}) = 0.
\end{equation}
\end{proof} 
\end{lemma} 

To simplify the computations in the subsequent analysis, we construct an adapted orthonormal frame. 
\begin{lemma}\label{lemma3.2}
Let 
$
\phi : (M^{n+1}(c), g) \longrightarrow (N^n, h)
$
be a biharmonic Riemannian submersion with an adapted frame 
$
\{e_1, e_2, \dots,e_n, e_{n+1}\}
$
and the integrability data 
$
\{f_{ij}^k, \kappa_i, \sigma_{ij} \} _{i,j,k = 1,\dots,n}.
$
Let $U$ be a nonempty open set on which $\kappa\neq0$. Then there exists
a nonempty open subset $V\subset U$ and an adapted orthonormal frame on
$V$ such that
\[
\kappa_2=\cdots=\kappa_n=0
\]
and
\[
\sigma_{i,i+j}=0 \qquad (i=1,\dots,n-2,\ j\geq2).
\]

\begin{proof}
We construct the desired frame locally. Throughout the proof, we shrink
the neighborhood $V$ whenever necessary.
Let 
\[
\kappa
=
\begin{pmatrix} 
\kappa_1\\
\kappa_2\\
\vdots\\
\kappa_n
\end{pmatrix},~
\sigma=(\sigma_{ij})_{n\times n,(1\leq i,j \leq n)},~
\text{and} ~
e^n_1=
\begin{pmatrix} 
1\\
0\\
\vdots\\
0
\end{pmatrix}(\in\mathbb{R}^n).
\]
We use Householder transformation. 
Set
\[
u^0=\frac{\kappa-\alpha_0||\kappa||e^n_1}{||\kappa-\alpha_0||\kappa||e^n_1||},
\]
where we take $\alpha_0\in\{-1,1\}$ such that $\kappa-\alpha_0||\kappa||e^n_1\not=0$.
Let 
\[
K^{(0)}
=I_n-2U^{(0)},
\]
where 
$U^{(0)}=(u^0_iu^0_j)_{n\times n,(1\leq i,j \leq n)}$.
Then, $K^{(0)}$ satisfies  
$
K^{(0)}\kappa=\alpha_0||\kappa||e_1^{n}
$
and 
$K^{(0)}\in O(n).$

Let
\[
S^{(0)}
=(s^0_{ij})
=K^{(0)}\sigma (K^{(0)})^T,~
x_{12}
=
\begin{pmatrix} 
s^0_{21}\\
s^0_{31}\\
\vdots\\
s^0_{n1}
\end{pmatrix},~
\text{and} ~
e^{n-1}_1=
\begin{pmatrix} 
1\\
0\\
\vdots\\
0
\end{pmatrix}(\in\mathbb{R}^{n-1}),
\]
where $(K^{(0)})^T$ denotes the transpose of $K^{(0)}$.
If $x_{12}=0$ on $V$, then we take $H^{(1)}=I_{n-1}$.
If $x_{12}$ is not identically zero on $V$, then, after replacing $V$
by a smaller nonempty open subset if necessary, we may assume that
$x_{12}\neq0$ on $V$. We then set
\[
u^1=\frac{x_{12}-\alpha_1||x_{12}||e^{n-1}_1}{||x_{12}-\alpha_1||x_{12}||e^{n-1}_1||},
\]
where we take $\alpha_1\in\{-1,1\}$ such that $x_{12}-\alpha_1||x_{12}||e^{n-1}_1\not=0$.
Let
\[
H^{(1)}
=
I_{n-1}-2U^{(1)},
\]
where 
$U^{(1)}=(u^1_iu^1_j)_{n-1 \times n-1,(1\leq i,j \leq n-1)}$.
Furthermore, let
\[
Q_1
=
\begin{pmatrix}
I_{1} & 0\\
0 & H^{(1)}
\end{pmatrix}~
\text{and}~
K^{(1)}
=Q_1K^{(0)}.~
\]
Obviously, $Q_1\in O(n)$.
We repeat the above argument.
For any $p=1,2,3,\cdots, n-3$,
let
\[
S^{(p)}
=(s^p_{ij})
=K^{(p)}\sigma (K^{(p)})^T,~
x_{p+1,p+2}
=
\begin{pmatrix} 
s^p_{p+2,p+1}\\
s^p_{p+3,p+1}\\
\vdots\\
s^p_{n,p+1}
\end{pmatrix},~
\text{and} ~
e^{n-(p+1)}_1=
\begin{pmatrix} 
1\\
0\\
\vdots\\
0
\end{pmatrix}(\in\mathbb{R}^{n-(p+1)}).
\]
If $x_{p+1,p+2}=0$ on $V$, then we take $H^{(p+1)}=I_{n-(p+1)}$.
If $x_{p+1,p+2}$ is not identically zero on $V$, then, after replacing
$V$ by a smaller nonempty open subset if necessary, we may assume that
$x_{p+1,p+2}\neq0$ on $V$. We then set
\[
u^{p+1}=\frac{x_{p+1,p+2}-\alpha_{p+1}||x_{p+1,p+2}||e^{n-(p+1)}_1}{||x_{p+1,p+2}-\alpha_{p+1}||x_{p+1,p+2}||e^{n-(p+1)}_1||},
\]
where we take $\alpha_{p+1}\in\{-1,1\}$ such that $x_{p+1,p+2}-\alpha_{p+1}||x_{p+1,p+2}||e^{n-(p+1)}_1\not=0$.
Let
\[
H^{(p+1)}
=
I_{n-(p+1)}-2U^{(p+1)},
\]
where 
$U^{(p+1)}=(u^{p+1}_iu^{p+1}_j)_{n-(p+1) \times n-(p+1),(1\leq i,j \leq n-(p+1))}$.
Furthermore, let
\[
Q_{p+1}
=
\begin{pmatrix}
I_{p+1} & 0\\
0 & H^{(p+1)}
\end{pmatrix},~
K^{(p+1)}
=Q_{p+1}K^{(p)},~
\text{and}~~
S^{(p+1)}
=Q_{p+1}S^{(p)}(Q_{p+1})^T.         
\]
Then 
$
\sigma'=(\sigma'_{ij})=K^{(n-2)}\sigma(K^{(n-2)})^T
$
satisfies that for any $i=1,2,3,\cdots ,n-2$, $\sigma'_{i,i+j}=0$ $(j\geq2)$.

One can also check that 
\[
K^{(n-2)}\kappa
=Q_{n-2}Q_{n-3}\cdots Q_{1}K^{(0)}\kappa
=
\begin{pmatrix} 
f\\
0\\
\vdots\\
0
\end{pmatrix},
\]
for some function $f$ and 
$K^{(n-2)}\in O(n)$, that is, 
\[
\begin{pmatrix}
e_1'\\
e_2'\\
\vdots\\
e_n'
\end{pmatrix}
=K^{(n-2)}
\begin{pmatrix}
e_1\\
e_2\\
\vdots\\
e_n
\end{pmatrix} 
\]
is the desired orthonormal frame.
Since the entries of $K^{(n-2)}=(K^{(n-2)}{}_{ij})$ are functions of the integrability data
$\kappa_i$ and $\sigma_{ij}$, Lemma~\ref{lemma3.1} implies that
$e_{n+1}(K^{(n-2)}{}_{ij})=0$.
Therefore the new orthonormal frame is adapted to the Riemannian submersion.
\end{proof}

\end{lemma}

The following lemma gives equations satisfied by a biharmonic Riemannian submersion in this frame.
\begin{lemma}\label{lemma3.3}
Let 
$
\phi : (M^{n+1}(c), g) \longrightarrow (N^n, h)
$
be a biharmonic Riemannian submersion with an adapted frame 
$
\{e_1, e_2,  \dots ,e_n,e_{n+1}\},
$
 the integrability data 
$
\{f_{ij}^k, \kappa_i, \sigma_{ij} \} _{i,j,k = 1,\dots,n}
$
with $\kappa_2 = \kappa_3 = \dots = \kappa_n =0$, and for any $i=1,\cdots ,n-2$, $\sigma_{i,i+j}=0,~(j\geq2)$. Then, the following equations hold.
\[
\left\{
\begin{array}{l}
\displaystyle
\Delta^M\kappa_1- \sum _ { i,l = 1 } ^ { n } \kappa _ { 1 } ( P_ { i 1 } ^ { l } ) ^ { 2 } +3\kappa_{1} \sigma_{12}^{2}+ ( n - 1 ) \kappa _ { 1 } c = 0 \quad (k=1), \\[2mm]

\displaystyle
\sum_{i=1}^{n} \left( 
2 e_i(\kappa_1) P_{i1}^k + 
\kappa_1 e_i(P_{i1}^k) + 
\kappa_1 P_{i1}^l P_{il}^k - 
\kappa_i \kappa_1 P_{i1}^k - 
\kappa_1 P_{ii}^l P_{l1}^k 
\right) 
-3\kappa_1\sigma_{12}\sigma_{2k}=0
\quad (k \neq 1).
\end{array}
\right.
\]
Here, the Laplacian is given by
\[
\Delta^M\kappa_1 = \sum_{i=1}^{n} e_i e_i(\kappa_1) - \sum_{i=1}^{n} \sum_{j=1}^{n} P_{ii}^j e_j(\kappa_1) - \kappa_1 e_1(\kappa_1).
\]

\begin{proof}
Using the assumptions on the adapted frame and substituting
$\kappa_2 = \kappa_3 = \dots = \kappa_n =0$ and for any $i=1,\cdots ,n-2$, $\sigma_{i,i+j}=0,~(j\geq2)$ into (\ref{eq:biha}), we obtain
\[
\Delta^M\kappa_1= \sum_{i=1}^{n} e_{i} e_{i}\left(\kappa_{1}\right)-\sum_{i, j=1}^{n} P_{i i}^{j} e_{j}\left(\kappa_{1}\right)-\kappa_{1} e_{1}\left(\kappa_{1}\right),  
\]
\[
\sum_{i,j=1}^n \left[ 2 e_i(\kappa_j) P_{ij}^1 + \kappa_j (e_i P_{ij}^1) + \kappa_j P_{ij}^l P_{il}^1 - \kappa_i \kappa_j P_{ij}^1 - \kappa_j P_{ii}^l P_{lj}^1 \right]= -\sum _ { i,l = 1 } ^ { n } \kappa _ { 1 } ( P_ { i 1 } ^ { l } ) ^ { 2 }, 
\]
\[
 \text{Ricci}^N(d\phi(\mu), d\phi(e_1))= 3\kappa_{1} \sigma_{12}^{2}+ ( n - 1 ) \kappa _ { 1 } c,
\]
and hence
\[
\Delta^M\kappa_1- \sum _ { i ,l = 1 } ^ { n } \kappa _ { 1 } ( P_ { i 1 } ^ { l } ) ^ { 2 } +3\kappa_{1} \sigma_{12}^{2}+ ( n - 1 ) \kappa _ { 1 } c = 0.
\]

 When $k \neq 1$, 

\[
\Delta^M\kappa_k= 0,
\]  
\begin{flalign*}
&\sum_{i,j=1}^n \left[ 2 e_i(\kappa_j) P_{ij}^k + \kappa_j (e_i P_{ij}^k) + \kappa_j P_{ij}^l P_{il}^k - \kappa_i \kappa_j P_{ij}^k - \kappa_j P_{ii}^l P_{lj}^k \right] &&\\
&= \sum_{i=1}^{n} \left( 2 e_i(\kappa_1) P_{i1}^k + \kappa_1 e_i(P_{i1}^k) + \kappa_1 P_{i1}^l P_{il}^k - \kappa_i \kappa_1 P_{i1}^k - \kappa_1 P_{ii}^l P_{l1}^k \right), &&
\end{flalign*}
\[
 \text{Ricci}^N(d\phi(\mu), d\phi(e_k))= -3\kappa_1\sigma_{12}\sigma_{2k},
\]
and hence
 \begin{align*}
 \sum_{i=1}^{n} \left( 2 e_i(\kappa_1) P_{i1}^k + \kappa_1 e_i(P_{i1}^k) + \kappa_1 P_{i1}^l P_{il}^k - \kappa_i \kappa_1 P_{i1}^k - \kappa_1 P_{ii}^l P_{l1}^k \right)
-3\kappa_1\sigma_{12}\sigma_{2k}=0.
 \end{align*}
\end{proof} 
\end{lemma}

We are now in a position to prove Theorem~\ref{main}.  
Lemmas~\ref{lemma3.1},~\ref{lemma3.2}, and ~\ref{lemma3.3} provide the necessary geometric and analytic tools for the argument.

\begin{theorem}[Theorem~\ref{main} (restated)]
For $n\geq2$, let 
$
\phi : (M^{n+1}(c), g) \longrightarrow (N^n, h)
$
be a Riemannian submersion from a Riemannian manifold with constant sectional curvature c. 
Then $ \phi $  is biharmonic if and only if it is harmonic.
\begin{proof} 
Since the case for $n=2$ has already been proven (cf. \cite{WO}), we will show the case for $n \geq 3$.

As stated in the definition and equation (33) in \cite{AO}, the tension field of the Riemannian submersion \(\phi\) is given by
\begin{align*}
\tau(\phi)& = \nabla ^\phi_{e_i} d\phi(e_i) - d\phi(\nabla^M_{e_i} e_i)\\
&= -d\phi(\nabla^M_{e_{n+1}} e_{n+1})\\
&= -d\phi((\nabla^M_{e_{n+1}} e_{n+1})^\mathcal{H})\\
&= -d\phi(\sum_{i=1}^n \kappa_i e_i)\\
&=- \sum_{i=1}^n \kappa_i \varepsilon_i.
\end{align*}

Since every harmonic map is biharmonic, it remains to prove the converse.
Assume that $\phi$ is biharmonic. We prove that $\phi$ is harmonic.
Suppose, to the contrary, that $\phi$ is not harmonic. Then there exists
a nonempty open set $\Omega_1$ on which $\kappa=(\kappa_1,\kappa_2,\cdots,\kappa_n)\neq0$. By Lemma~\ref{lemma3.2},
after replacing $\Omega_1$ by a smaller nonempty open subset if necessary, we may
choose an adapted orthonormal frame \(\{e_1, e_2,\dots, e_n,e_{n+1}\}\) such that
$
\kappa_2=\cdots=\kappa_n=0
$
and
$
\sigma_{i,i+j}=0\ (i=1,\dots,n-2,\ j\geq2).
$
Consequently, the tension field of the Riemannian submersion \(\phi\) can be expressed as
\[
\tau(\phi) = - \kappa_1 \varepsilon_1.
\]

In terms of this frame, the curvature equation (\ref{eq:R^M}) simplifies to the following expression:
\begin{equation}
\left\{ \,
\begin{aligned}
& R_{i(n+1) 1(n+1)}^{M} :\quad e_{i}\left(\kappa_{1}\right)=-\sigma_{12} \sigma_{i 2}+\delta_{1 i}\left(\kappa_{1}^{2}+c\right),\\
& R_{1(n+1) 1(n+1)}^{M} :\quad e_{1}\left(\kappa_{1}\right)=-\sigma_{12}^{2}+\kappa_{1}^{2}+c,\\
& R_{i(n+1) 12}^{M} :\quad e_{i}\left(\sigma_{12}\right)=P_{i 1}^{l} \sigma_{l 2}+\kappa_{1} \sigma_{i 2}+\delta_{i 1} \kappa_{1} \sigma_{12},\\
&R_{i(n+1) m(n+1)}^{M} :\quad P_{i m}^{1}=\frac{1}{\kappa_{1}}\left(\sigma_{ml}  \sigma_{i l}-\delta_{m i} c\right)\quad(m \neq 1),\\
& R_{i(n+1) j 1}^M :\quad P_{i j}^{2} \sigma_{12}=P_{i l}^{1} \sigma_{l j}-\kappa_{1} \sigma_{i j}\quad(j \geq 3),\\
& R_{1(n+1) 2 m}^{M} :\quad e_{1}\left(\sigma_{2 m}\right)=P_{12}^{l} \sigma_{l m}-P_{1 m}^{l} \sigma_{l 2}-\kappa_{m} \sigma_{12}+\kappa_{1} \sigma_{2 m},\\
&R _ { 1 (n + 1)im}^M :\qquad e _ { 1 } ( \sigma _ { im} ) =-P_ { 1m} ^ { l}  \sigma _ { li}+P_ { 1i } ^ { l} \sigma _ { lm } -\kappa _ { m}\sigma _ { 1i } + \kappa _ { 1 } \sigma _ { im} +\kappa _ { i} \sigma _ { 1m },\\ 
&R _ { 1 (n + 1)im}^M :\qquad e _ { 1 } ( \sigma _ { im} ) =-P_ { 1m} ^ { l}  \sigma _ { li}+P_ { 1i } ^ { l} \sigma _ { lm } + \kappa _ { 1 } \sigma _ { im} \quad(i,m\ge3).
\end{aligned}
\label{eq:R^M3}
\right.
\end{equation} 
Consequently, we have
\begin{equation}
\begin{aligned}
\sum_{i,m=1}^n e _ { 1 } ( \sigma _ { i m } ^ { 2 } ) &=4 \kappa_{1} \sigma_{12}^{2}+2 \kappa_{1}\sum_{i,m=1}^n \sigma_{i m}^{2}, \\
 e_1(\sigma _{12}^2)
&=2\frac{\sigma_{12}^2}{\kappa_1}\Big(2\kappa_1^2-\sum_{l=2}^n\sigma_{l2}^2\Big).
\end{aligned}
\end{equation}

By Lemma \ref{lemma3.3}, the first term of the biharmonic equation satisfies
\begin{align*}
& \sum_{i=1}^{n} e_{i} e_{i}\left(\kappa_{1}\right)-\sum_{i, j=1}^{n} P_{i i}^{j} e_{j}\left(\kappa_{1}\right)-\kappa_{1} e_{1}\left(\kappa_{1}\right)\\
 &\qquad\qquad\qquad- \sum _ { i ,l= 1 } ^ { n } \kappa _ { 1 } ( P_ { i 1 } ^ { l } ) ^ { 2 } +3\kappa_{1} \sigma_{12}^{2}+ ( n - 1 ) \kappa _ { 1 } c = 0.
\end{align*}
Substituting equation (\ref{eq:R^M3}) into the remaining terms yields the following identities:
\begin{align*}
 & \sum_{i=1}^{n} e_{i} e_{i}(\kappa_{1}) -\kappa_{1} e_{1}(\kappa_{1}) \\
 &\quad = \kappa_{1} \left( e_{1}(\kappa_{1}) - 3 \sigma_{12}^{2} - \sum_{i=1}^{n} \sigma_{2i}^{2} \right)  \\
 &\qquad\qquad+ \sum_{i,l=1}^{n} \left( P_{ii}^{l} \sigma_{12} \sigma_{2l} + P_{i2}^{l} \sigma_{12} \sigma_{li} + P_{i1}^{l} \sigma_{2i} \sigma_{l2} \right) \\
&\quad= \kappa_{1} \left( - 4 \sigma_{12}^{2}+\kappa_1^2+c - \sum_{i=1}^{n} \sigma_{2i}^{2} \right)+ \sum_{i,l=1}^{n} P_{ii}^{l} \sigma_{12} \sigma_{2l} \\
&\qquad - \frac{1}{\kappa_{1}} \left( \sum_{i,p=1}^{n}\sum_{m=3}^{n}  \sum_{l=2}^{n} \sigma_{lp} \sigma_{ip} \sigma_{lm} \sigma_{mi}\right. \\
&\qquad\qquad\qquad\qquad\qquad\qquad
\left.+\sum_{i=2}^{n}\sum_{m=3}^{n}\sigma_{mi}^{2} c -\sum_{l=1}^n \sigma_{2l}^2\sigma_{12}^2+\sigma_{12}^2c \right)\\
&\qquad - \sum_{i=1}^{n}\sum_{m=3}^{n} \kappa_{1} \sigma_{im}^{2} 
+ \frac{1}{\kappa_{1}} \left(- \sum_{i,p=1}^{n} \sum_{l=2}^{n}  \sigma_{lp} \sigma_{ip} \sigma_{2i} \sigma_{l2} - \sum_{l=2}^{n} \sigma_{2l}^{2} c \right),
\end{align*}
\begin{align*}
&-\sum_{i, j=1}^{n} P_{i i}^{j} e_{j}\left(\kappa_{1}\right)
=\sum_{i, l=1}^{n} P_{i i}^l \sigma_{12} \sigma_{l 2}-\sum_{i=2}^{n}\sum_{l=1}^{n} \kappa_{1} \sigma_{i l}^{2} \\
&\qquad\qquad\qquad\qquad
+(n-1) \kappa_{1} c +\frac{1}{\kappa_{1}}\left\{-\sum_{i=2}^{n}\sum_{l=1}^{n}\sigma_{i l}^{2} c+(n-1) c^{2}\right\},
\end{align*}
\begin{align*}
& - \sum _ { i,l = 1 } ^ { n } \kappa _ { 1 } ( P_ { i 1 } ^ { l } ) ^ { 2 }  =-\sum_{i=1}^{n} \sum_{l=2}^{n} \frac{1}{\kappa_{1}}\left(\sum_{p=1}^{n} \sigma_{l p} \sigma_{i p}-\delta_{l i} c\right)^{2} \\ 
&\qquad \qquad \qquad 
=\frac{1}{\kappa_{1}}\left(-\sum_{i,p,q=1}^{n} \sum_{l=2}^{n}  \sigma_{l p} \sigma_{i p} \sigma_{l q} \sigma_{i q}+\sum_{i=2}^{n} \sum_{p=1}^{n} 2 \sigma_{i p}^{2} c-(n-1) c^{2}\right).
\end{align*}
From these relations, we deduce that\\
\[
\kappa_1 ( \kappa_{1}^{2} +\sigma_{12}^2-2\sum_{i, m=1}^{n} \sigma_{i m}^{2}+(2 n-1)c )=0.
\]
Since \( \kappa_1 \neq 0 \), it follows that
\begin{equation}\label{key}
\kappa_{1}^{2} +\sigma_{12}^2-2\sum_{i, m=1}^{n} \sigma_{i m}^{2}+(2 n-1)c=0.
\end{equation}

By differentiating both sides of \eqref{key} by $e_1$ and substituting equation \eqref{key} again, we have
\begin{equation}\label{e1}
4\sigma_{12}^2+\frac{\sigma_{12}^2}{\kappa_1^2}\sum_{l=2}^{n} \sigma_{l2}^{2}+(2 n-2)c=0.
\end{equation}

If $c>0$, then we have a contradiction.

Assume that $c=0$, then we have $\sigma_{12}\equiv0$. 
For $j=2,3,\cdots ,n-1$, by the first equation of \eqref{eq:R^M}, 
\[
0=e_{j}(\sigma_{1,j+1})=P_{j1}^l\sigma_{l,j+1}+\kappa_1\sigma_{j,j+1}.
\]
By the fourth equation of \eqref{eq:R^M3}, we have
\[
\sigma_{j,j+1}
\left(
-\sigma_{j-1,j}^2-\sigma_{j,j+1}^2-\sigma_{j+1,j+2}^2+\kappa_1^2
\right)
=0~~\text{for}~n\geq4~\text{and}~j=2,3,4,\cdots, n-2,
\]
and 
\[
\sigma_{n-1,n}
\left(
-\sigma_{n-2,n-1}^2-\sigma_{n-1,n}^2+\kappa_1^2
\right)
=0.
\]

Assume that $\sigma_{23}\not=0$ at some point $p_{2}$, that is $\sigma_{23}\not=0$ on some open set $\Omega_{2}\subset \Omega_1$.
Then, one has 
\[
\kappa_1^2=\sigma_{23}^2+\sigma_{34}^2.
\]
Substituting this into \eqref{key}, we have a contradiction. Hence one has $\sigma_{23}=0$.

We iterate the same argument. 
Assume that $\sigma_{j,j+1}\not=0$ at some point $p_{j}$, that is $\sigma_{j,j+1}\not=0$ on some open set $\Omega_{j}$.
Then, one has 
\[
\kappa_1^2=\sigma_{j,j+1}^2+\sigma_{j+1,j+2}^2.
\]
Substituting this into \eqref{key}, we have a contradiction. Hence one has $\sigma_{j,j+1}=0$.
 Finally we have
 \[
 \sigma_{n-1,n}
 \left(
-\sigma_{n-1,n}^2+\kappa_1^2
 \right)=0.
 \]
 
Assume that $\sigma_{n-1,n}\not=0$ at some point $p_{n-1}$, that is $\sigma_{n-1,n}\not=0$ on some open set $\Omega_{n-1}$.
Then, one has 
\[
\kappa_1^2=\sigma_{n-1,n}^2.
\]
Substituting this into \eqref{key}, we have a contradiction.
Hence one has $\sigma_{n-1,n}=0$.
Therefore, $\sigma_{ij}=0$ for all $i,j$. By \eqref{key}, we have $\kappa_{1}=0$, which is a contradiction.

Therefore, we have $c<0.$
When $n=3$, we use the convention $\sigma_{34}=0$.
Thus the following argument also applies to the case $n=3$.

By differentiating both sides of \eqref{key} by $e_3$, we have
\[
\sigma_{12}\sigma_{23}(\sigma_{23}^2+\sigma_{34}^2-c+4\kappa_1^2)=0.
\]
Since
\[
\sigma_{23}^2+\sigma_{34}^2-c+4\kappa_1^2>0,
\]
we have 
\[
\sigma_{12}\sigma_{23}=0.
\]

By \eqref{e1}, $\sigma_{12}$ is a non-zero constant.
By differentiating $\sigma_{12}$ by $e_1$, we have $\kappa_1\sigma_{12}=0$, which is a contradiction.

\end{proof} 
\end{theorem}



\noindent
{\bf Acknowledgements.}~\\
The first author is partially supported by the Grant-in-Aid for Scientific Research (C), No.23K03107, Japan Society for the Promotion of Science.




\bibliographystyle{amsbook}

\end{document}